\documentclass{ifacconf}

\makeatletter
\let\old@ssect\@ssect 
\makeatother

\usepackage{graphicx}      
\usepackage{natbib}        
\usepackage{amsmath,amsfonts,amssymb}
\usepackage{lmodern}
\usepackage{hyperref}
\usepackage{tikz}
\usepackage{color}
\usepackage{algorithm}
\usepackage{algorithmicx}
\usepackage{algpseudocode}
\algnewcommand\algorithmicinput{\textbf{Input:}}
\algnewcommand\Input{\item[\algorithmicinput]}
\algnewcommand\algorithmicoutput{\textbf{Output:}}
\algnewcommand\Output{\item[\algorithmicoutput]}

\newcommand{\N}{\mathbb{N}}
\newcommand{\R}{\mathbb{R}}

\newcommand{\norm}[1]{\Vert #1 \Vert}
\DeclareMathOperator*{\argmin}{argmin}
\newcommand{\s}{\mathrm{s}}

\newcommand{\proj}[2]{P_{#1}(#2)}

\newcommand{\gencone}[4]{{#1}_{#2}^{#4}(#3)} 
\newcommand{\tancone}[2]{\gencone{T}{#1}{#2}{}} 

\newcommand{\ushort}[1]{\mkern 0.9mu\underline{\mkern-0.9mu#1\mkern-0.9mu}\mkern 0.9mu}

\DeclareMathOperator{\rank}{rank}
\DeclareMathOperator{\diag}{diag}

\newcommand{\pgd}{\mathrm{P}^2\mathrm{GD}}
\newcommand{\pgdr}{\mathrm{P}^2\mathrm{GDR}}

\makeatletter
\def\@ssect#1#2#3#4#5#6{%
  \NR@gettitle{#6}
  \old@ssect{#1}{#2}{#3}{#4}{#5}{#6}
}
\makeatother

\edef\endfrontmatter{%
  \unexpanded\expandafter{\endfrontmatter}
  \noexpand\endNoHyper 
}

\begin{document}
\begin{frontmatter}

\title{Comparison of an Apocalypse-Free and an Apocalypse-Prone First-Order Low-Rank Optimization Algorithm\thanksref{footnoteinfo}} 

\thanks[footnoteinfo]{This work was supported by the Fonds de la Recherche Scientifique -- FNRS and the Fonds Wetenschappelijk Onderzoek -- Vlaanderen under EOS Project no 30468160. K. A. Gallivan is partially supported by the U.S. National Science Foundation under grant CIBR 1934157. This work has been submitted to IFAC for possible publication.}

\author[First]{Guillaume Olikier} 
\author[Second]{Kyle A. Gallivan} 
\author[First]{P.-A. Absil}

\address[First]{ICTEAM Institute, UCLouvain, Avenue Georges Lema\^{\i}tre 4, 1348 Louvain-la-Neuve, Belgium (\href{mailto:guillaume.olikier@uclouvain.be}{\nolinkurl{guillaume.olikier@uclouvain.be}}, \href{mailto:pa.absil@uclouvain.be}{\nolinkurl{pa.absil@uclouvain.be}}).}
\address[Second]{Department of Mathematics, Florida State University, 1017 Academic Way, Tallahassee, FL 32306-4510, USA (\href{mailto:kgallivan@fsu.edu}{\nolinkurl{kgallivan@fsu.edu}}).}

\begin{abstract}                
We compare two first-order low-rank optimization algorithms, namely $\pgd$ (Schneider and Uschmajew, 2015), which has been proven to be apocalypse-prone (Levin et al., 2021), and its apocalypse-free version $\pgdr$ obtained by equipping $\pgd$ with a suitable rank reduction mechanism (Olikier et al., 2022). Here an apocalypse refers to the situation where the stationarity measure goes to zero along a convergent sequence whereas it is nonzero at the limit.
The comparison is conducted on two simple examples of apocalypses, the original one (Levin et al., 2021) and a new one. We also present a potential side effect of the rank reduction mechanism of $\pgdr$ and discuss the choice of the rank reduction parameter.
\end{abstract}


\end{frontmatter}

\section{Introduction}
As in \cite{OlikierGallivanAbsil2022}, we consider the problem
\begin{equation}
\label{eq:LowRankOpti}
\min_{X \in \R_{\le r}^{m \times n}} f(X)
\end{equation}
of minimizing a differentiable function $f : \R^{m \times n} \to \R$ with locally Lipschitz continuous gradient on the \emph{determinantal variety} \cite[Lecture~9]{Harris}
\begin{equation*}
\R_{\le r}^{m \times n} := \{X \in \R^{m \times n} \mid \rank X \le r\},
\end{equation*}
$m$, $n$, and $r$ being positive integers such that $r < \min\{m,n\}$.
This problem appears in several applications such as matrix equations, model reduction, matrix sensing, and matrix completion; see, e.g., \cite{SchneiderUschmajew2015}, \cite{HaLiuFoygel2020}, and the references therein.
As problem \eqref{eq:LowRankOpti} is in general intractable---see \cite{GillisGlineur2011}---, our goal is to find a \emph{stationary point} of this problem, i.e., a zero of the \emph{stationarity measure}
\begin{equation*}
\s_f :
\R_{\le r}^{m \times n} \to \R :
X \mapsto \norm{\proj{\tancone{\R_{\le r}^{m \times n}}{X}}{-\nabla f(X)}},
\end{equation*}
denoted by $g^-$ in \cite{SchneiderUschmajew2015}, that returns the norm of any projection of $-\nabla f(X)$ onto the tangent cone to $\R_{\le r}^{m \times n}$ at $X$; the notation is introduced in Section~\ref{sec:NotationPreliminaries}.

To the best of our knowledge, the second-order method given in \cite[Algorithm~3.1]{LevinKileelBoumal2021} and the first-order method given in \cite[Algorithm~2]{OlikierGallivanAbsil2022} are the only two algorithms in the literature that provably converge to stationary points. Other algorithms, such as \cite[Algorithm~3]{SchneiderUschmajew2015}, can fail in the sense that they can produce a feasible sequence $(X_i)_{i \in \N}$ that converges to some point $X$ with the property that $\lim_{i \to \infty} \s_f(X_i) = 0 < \s_f(X)$. Such a triplet $(X, (X_i)_{i \in \N}, f)$ is called an \emph{apocalypse} and the point $X$, which necessarily satisfies $\rank X < r$, is said to be \emph{apocalyptic} according to \cite[Definition~2.7]{LevinKileelBoumal2021}.

In this paper, using synthetic instances of \eqref{eq:LowRankOpti}, we compare the behavior of the first-order algorithms given in \cite[Algorithm~3]{SchneiderUschmajew2015} and \cite[Algorithm~2]{OlikierGallivanAbsil2022}, respectively dubbed $\pgd$ and $\pgdr$, the latter consisting of the former equipped with a suitable rank reduction mechanism. We observe that the experiments corroborate the theory but also reveal that the choice of the rank reduction parameter can be significant in practice; see Section~\ref{sec:Conclusion} for details.

This paper is organized as follows. After recalling some notation and preliminaries in Section~\ref{sec:NotationPreliminaries}, we compare in Section~\ref{sec:TwoExamplesApocalypses} the two algorithms on two simple examples of apocalypses, one on $\R_{\le 2}^{3 \times 3}$ proposed in \cite[\S 2.2]{LevinKileelBoumal2021} and one on $\R_{\le 1}^{2 \times 2}$, to illustrate how $\pgdr$ avoids following the apocalypses due to its rank reduction mechanism. In Section~\ref{sec:PotentialSideEffectRankReductionMechanism}, we present a potential side effect of that mechanism. We discuss the choice of the rank reduction parameter and draw conclusions in Section~\ref{sec:Conclusion}.

\section{Notation and preliminaries}
\label{sec:NotationPreliminaries}
In this section, we recall some notation and preliminaries from \cite{OlikierGallivanAbsil2022} to which we refer for a more complete review of the background material.
In what follows, $\R^{m \times n}$ is endowed with the Frobenius inner product, and $\norm{\cdot}$ denotes the Frobenius norm.
A nonempty subset $\mathcal{C}$ of $\R^{m \times n}$ is said to be a \emph{cone} if, for every $X \in \mathcal{C}$ and every $\lambda \in [0,\infty)$, it holds that $\lambda X \in \mathcal{C}$.
For every nonempty subset $\mathcal{S}$ of $\R^{m \times n}$ and every $X \in \mathcal{S}$, the set $\tancone{\mathcal{S}}{X}$ of all $V \in \R^{m \times n}$ such that there exist $(t_i)_{i \in \N}$ in $(0,\infty)$ converging to $0$ and $(V_i)_{i \in \N}$ in $\R^{m \times n}$ converging to $V$ such that $X+t_iV_i \in \mathcal{S}$ for every $i \in \N$ is a closed cone, not necessarily convex however, called the \emph{tangent cone} to $\mathcal{S}$ at $X$.
For every closed cone $\mathcal{C}$ in $\R^{m \times n}$ and every $X \in \R^{m \times n}$, the set $\proj{\mathcal{C}}{X} := \argmin_{Y \in \mathcal{C}} \norm{X-Y}$, called the \emph{projection} of $X$ onto $\mathcal{C}$, is nonempty, compact, and all its elements have the same norm. If $\argmin$ is a singleton, we identify it with its element.

The iteration map of $\pgd$ \cite[Algorithm~1]{OlikierGallivanAbsil2022} is given as Algorithm~\ref{algo:P2GDmap}; the only difference with the one of \cite[Algorithm~3]{SchneiderUschmajew2015} is that the initial step size for the backtracking procedure is chosen in a given bounded interval and not in $[1,\infty)$. The acronym ``$\pgd$'' follows from the fact that the iteration map of this algorithm consists of a step along a projection of the negative gradient onto the tangent cone to $\R_{\le r}^{m \times n}$, followed by a projection onto $\R_{\le r}^{m \times n}$. The first projection can be computed by \cite[Algorithm~2]{SchneiderUschmajew2015} and the second can be obtained by truncating an SVD in view of the Eckart--Young theorem. Because of the ``Choose'' statements, the $\pgd$ map is set-valued in general. In what follows, $\pgd(X; f, \ushort{\alpha}, \bar{\alpha}, \beta, c)$ denotes the set of all possible outputs of Algorithm~\ref{algo:P2GDmap}.

\begin{algorithm}[H]
\caption{$\pgd$ map}
\label{algo:P2GDmap}
\begin{algorithmic}[1]
\Require
$(f, \ushort{\alpha}, \bar{\alpha}, \beta, c)$ where $f : \R^{m \times n} \to \R$ is differentiable with locally Lipschitz continuous gradient, $0 < \ushort{\alpha} \le \bar{\alpha} < \infty$, and $\beta, c \in (0,1)$.
\Input
$X \in \R_{\le r}^{m \times n}$ such that $\s_f(X) > 0$.
\Output
$Y \in \pgd(X; f, \ushort{\alpha}, \bar{\alpha}, \beta, c)$.

\State
Choose $G \in \proj{\tancone{\R_{\le r}^{m \times n}}{X}}{-\nabla f(X)}$, $\alpha \in [\ushort{\alpha},\bar{\alpha}]$, and $Y \in \proj{\R_{\le r}^{m \times n}}{X + \alpha G}$;
\While
{$f(Y) > f(X) - c \, \alpha \, \s_f(X)^2$}
\State
$\alpha \gets \alpha \beta$;
\State
Choose $Y \in \proj{\R_{\le r}^{m \times n}}{X + \alpha G}$;
\EndWhile
\State
Return $Y$.
\end{algorithmic}
\end{algorithm}

As mentioned above, $\pgdr$ consists of $\pgd$ equipped with a rank reduction mechanism, hence the ``R'' in the acronym. This mechanism uses the numerical rank: given $\Delta \in [0,\infty)$ and $X \in \R^{m \times n} \setminus \{0_{m \times n}\}$, the \emph{$\Delta$-rank} of $X$ is defined as
\begin{equation}
\label{eq:DeltaRank}
\rank_\Delta X := \max\{j \in \{1, \dots, \rank X\} \mid \sigma_j(X) > \Delta\},
\end{equation}
where $\sigma_1(X) \ge \dots \ge \sigma_{\min\{m,n\}}(X)$ denote the singular values of $X$, and the definition is completed by setting $\rank_\Delta 0_{m \times n} := 0$.
Based on this definition, the iteration map of $\pgdr$ \cite[Algorithm~3]{OlikierGallivanAbsil2022} is given as Algorithm~\ref{algo:P2GDRmap}.
In particular, $\pgd$ corresponds to $\pgdr$ with $\Delta := 0$, and, more generally, the smaller $\Delta$ is, the more $\pgdr$ tends to behave as $\pgd$. Furthermore, by the Eckart--Young theorem, for every $X \in \R^{m \times n}$ and every $\ushort{r} \in \{0, \dots, \rank X\}$,
\begin{equation*}
\proj{\R_{\ushort{r}}^{m \times n}}{X} = \proj{\R_{\le \ushort{r}}^{m \times n}}{X}.
\end{equation*}
As the $\pgd$ map, the $\pgdr$ map is set-valued in general. In what follows, $\pgdr(X; f, \ushort{\alpha}, \bar{\alpha}, \beta, c, \Delta)$ denotes the set of all possible outputs of Algorithm~\ref{algo:P2GDRmap}.

In all discussions and experiments, we use a constant initial step size for the backtracking line search, which amounts to choosing $\bar{\alpha} = \ushort{\alpha}$.

\begin{algorithm}[H]
\caption{$\pgdr$ map}
\label{algo:P2GDRmap}
\begin{algorithmic}[1]
\Require
$(f, \ushort{\alpha}, \bar{\alpha}, \beta, c, \Delta)$ where $f : \R^{m \times n} \to \R$ is differentiable with locally Lipschitz continuous gradient, $0 < \ushort{\alpha} \le \bar{\alpha} < \infty$, $\beta, c \in (0,1)$, and $\Delta \in (0,\infty)$.
\Input
$X \in \R_{\le r}^{m \times n}$ such that $\s_f(X) > 0$.
\Output
$Y \in \pgdr(X; f, \ushort{\alpha}, \bar{\alpha}, \beta, c, \Delta)$.

\For
{$j \in \{0, \dots, \rank X - \rank_\Delta X\}$}
\State
Choose $\hat{X}^j \in \proj{\R_{\rank X - j}^{m \times n}}{X}$;
\State
Choose $\tilde{X}^j \in \hyperref[algo:P2GDmap]{\pgd}(\hat{X}^j; f, \ushort{\alpha}, \bar{\alpha}, \beta, c)$;
\EndFor
\State
Return $Y \in \argmin_{\{\tilde{X}^j \mid j \in \{0, \dots, \rank X - \rank_\Delta X\}\}} f$.
\end{algorithmic}
\end{algorithm}

The iterative process is summarized in Algorithm~\ref{algo:IterativeP2GDR}, where $\Delta := 0$ corresponds to the $\pgd$ algorithm and $\Delta > 0$ to the $\pgdr$ algorithm. Although the convergence analysis in \cite{OlikierGallivanAbsil2022} is conducted for $\varepsilon := 0$, it is necessary to choose $\varepsilon > 0$ in a practical implementation to guarantee that the algorithm terminates after a finite number of iterations.

By \cite[Theorem~5.2 and Corollary~5.3]{OlikierGallivanAbsil2022}, if $\Delta > 0$ and $\varepsilon := 0$, then Algorithm~\ref{algo:IterativeP2GDR} produces either a finite sequence the last term of which is stationary or an infinite sequence with the following two properties: its accumulation points are stationary and the stationarity measure $\s_f$ goes to zero along each convergent subsequence. Thus, in the second case, except if the sequence diverges to infinity, choosing $\varepsilon > 0$ makes Algorithm~\ref{algo:IterativeP2GDR} terminate after finitely many iterations.
We do not have this guarantee if $\Delta := 0$. Indeed, if $\varepsilon := \Delta := 0$ and Algorithm~\ref{algo:IterativeP2GDR} produces a sequence $(X_i)_{i \in \N}$ that does not diverge to infinity, then, to the best of our knowledge, it is not known whether $\liminf_{i \to \infty} \s_f(X_i) = 0$.

\begin{algorithm}[H]
\caption{Iterative $\mathrm{P}^2\mathrm{GD(R)}$}
\label{algo:IterativeP2GDR}
\begin{algorithmic}[1]
\Require
$(X_0, f, \ushort{\alpha}, \bar{\alpha}, \beta, c, \Delta, \varepsilon)$ where $X_0 \in \R_{\le r}^{m \times n}$, $f : \R^{m \times n} \to \R$ is differentiable with locally Lipschitz continuous gradient, $0 < \ushort{\alpha} \le \bar{\alpha} < \infty$, $\beta, c \in (0,1)$, and $\Delta, \varepsilon \in [0,\infty)$.

\State
$i \gets 0$;
\While
{$\s_f(X_i) > \varepsilon$}
\State
Choose $X_{i+1} \in \hyperref[algo:P2GDRmap]{\pgdr}(X_i; f, \ushort{\alpha}, \bar{\alpha}, \beta, c, \Delta)$;
\State
$i \gets i+1$;
\EndWhile
\end{algorithmic}
\end{algorithm}

\section{Two examples of apocalypses}
\label{sec:TwoExamplesApocalypses}
In this section, we compare the behavior of $\pgd$ and $\pgdr$ on two examples of apocalypses. In Section~\ref{subsec:ExampleApocalypseLevinEtAl}, we compare the two algorithms empirically on the example of \cite[\S 2.2]{LevinKileelBoumal2021}. In Section~\ref{subsec:ExampleApocalypseSmallestSize}, we compare them analytically on a simple example of an apocalypse on $\R_{\le 1}^{2 \times 2}$.

\subsection{The example of Levin et al.}
\label{subsec:ExampleApocalypseLevinEtAl}
In \cite[\S 2.2]{LevinKileelBoumal2021}, the following instance of \eqref{eq:LowRankOpti} is considered: minimizing
\begin{equation*}
f : \R^{3 \times 3} \to \R : X \mapsto Q(X_{1:2,1:2}) + \phi(X_{3,3})
\end{equation*}
on $\R_{\le 2}^{3 \times 3}$, where $X_{1:2,1:2}$ is the upper-left $2 \times 2$ submatrix of $X$, $X_{3,3}$ its bottom-right entry, $\phi : \R \to \R : x \mapsto \frac{x^4}{4}-\frac{(x+1)^2}{2}$, $Q : \R^{2 \times 2} \to \R : Y \mapsto \frac{1}{2}\norm{D(Y-Y^*)}^2$, $D := \diag(1,\frac{1}{2})$, and $Y^* := \diag(1,0)$.
First, it is observed that $\argmin f = \diag(1,0,x_0) =: X^*$, where $x_0 := \argmin \phi \approx 1.32471795724475$, and $f^* := \min f = f(X^*) = \phi(x_0) \approx -1.932257884495233$.
Second, it is proven analytically that $\pgd$ follows an apocalypse if used on this problem with $X_0 := \diag(2,1,0)$, $\bar{\alpha} := \ushort{\alpha} := \frac{8}{5}$, any $\beta \in (0,1)$, and $c := \frac{1}{5}$.

\begin{figure}[h]
\begin{center}
\includegraphics[scale=0.6]{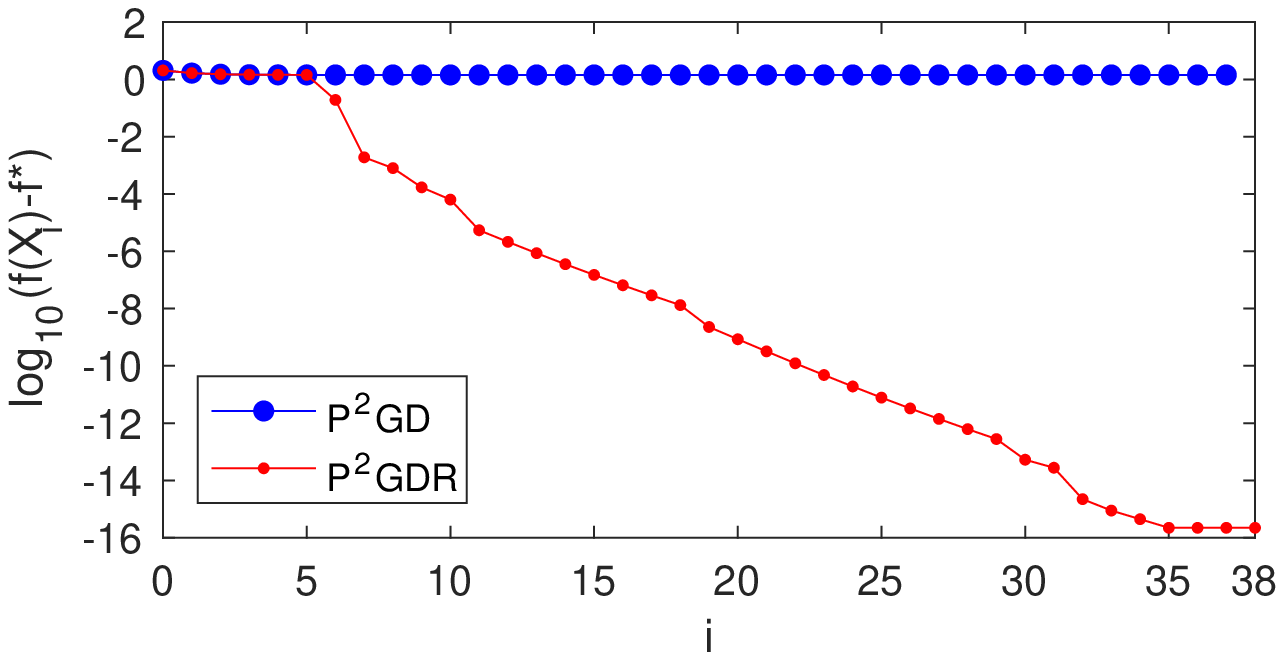}
\includegraphics[scale=0.6]{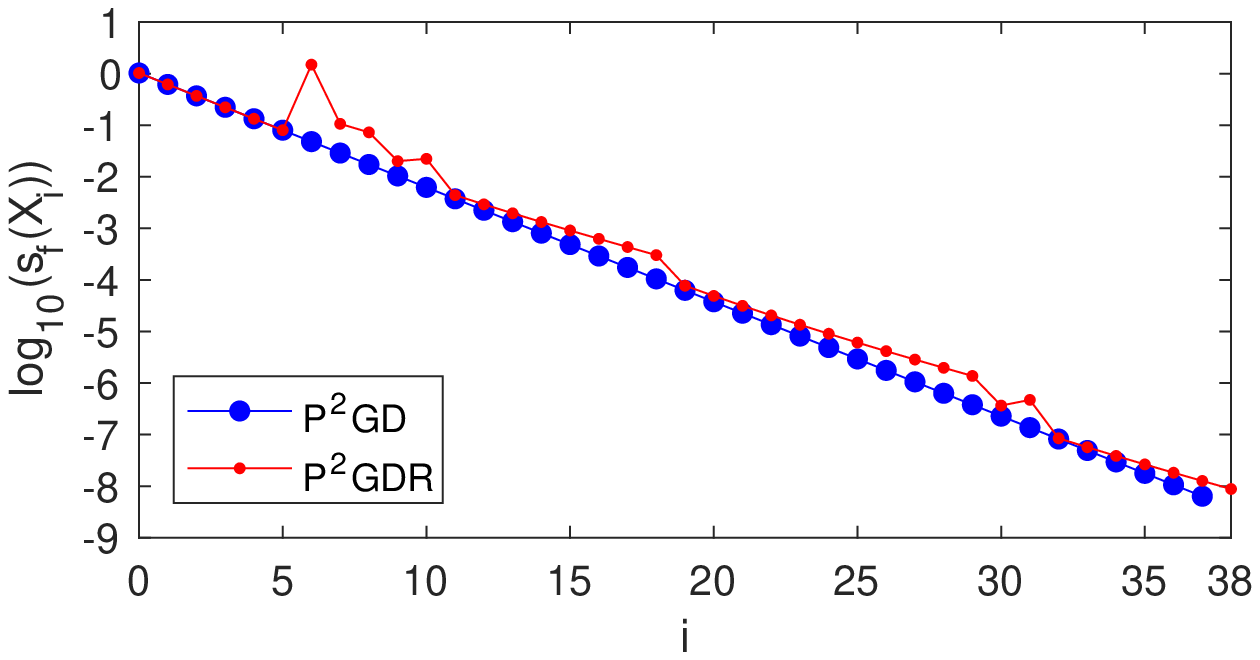}
\hphantom{\includegraphics[scale=0.6]{LevinExample_Cost.eps}}
\includegraphics[scale=0.6]{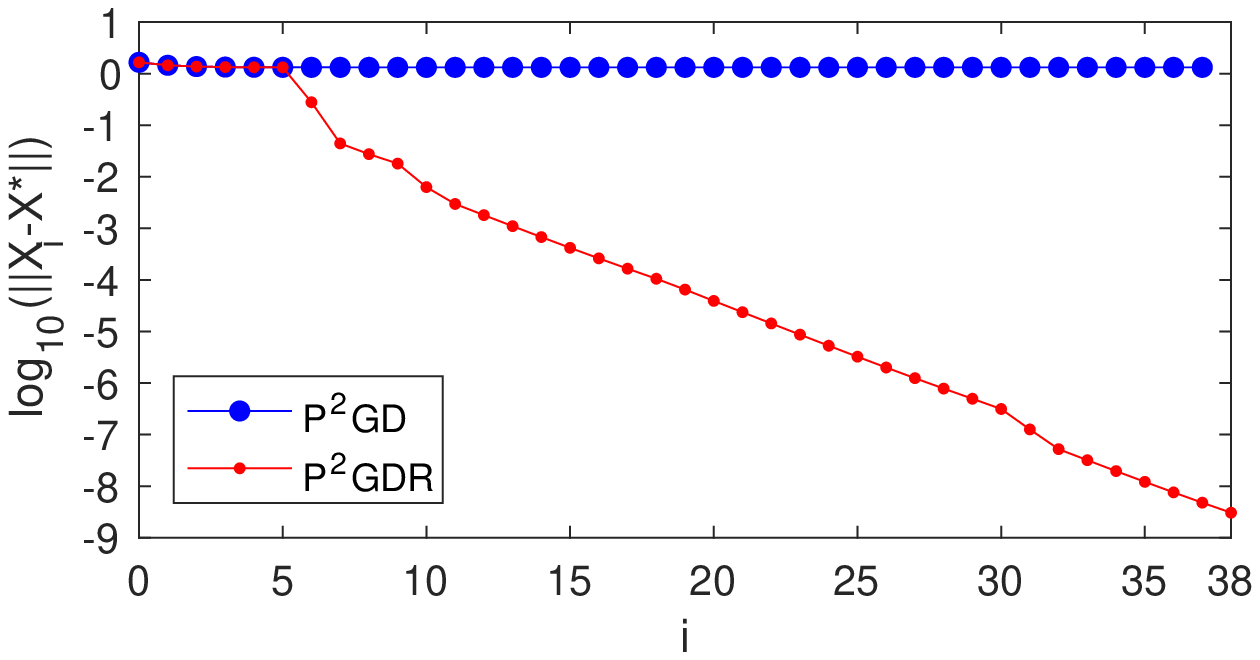}
\end{center}
\caption{Sequences $(f(X_i)-f^*)_{i=0}^{i=38}$, $(\s_f(X_i))_{i=0}^{i=38}$, and $(\norm{X_i-X^*})_{i=0}^{i=38}$ for the problem of Section~\ref{subsec:ExampleApocalypseLevinEtAl}.}
\label{fig:RunP2GDR}
\end{figure}

In this subsection, we verify empirically that, on the same problem with the same input parameters, $\pgd$ and $\pgdr$ with $\Delta := 0.1$ respectively behave as predicted in \cite[\S 2.2]{LevinKileelBoumal2021} and in agreement with the following theoretical guarantee: by \cite[Theorem~5.2 and Corollary~5.3]{OlikierGallivanAbsil2022}, since all sublevel sets of $f$ are bounded, for every initial iterate, $\pgdr$ produces a bounded sequence in $\R_{\le 2}^{3 \times 3}$ the accumulation points of which are stationary and along which the stationarity measure $\s_f$ goes to zero.
We run our Matlab implementation\footnote{Available at \url{https://sites.uclouvain.be/absil/2022.02}.} of Algorithm~\ref{algo:IterativeP2GDR} with $\beta := \frac{1}{2}$ and $\varepsilon := 10^{-8}$. This gives us the sequence $(X_i)_{i=0}^{i=37}$ for $\pgd$, where
\begin{equation}
\label{eq:ExampleApocalypseLevinEtAl}
X_i := \diag(1+(-3/5)^i,(3/5)^i,0)
\end{equation}
for every $i \in \{0, \dots, 37\}$. The sequence $(X_i)_{i=0}^{i=38}$ produced by $\pgdr$ obeys \eqref{eq:ExampleApocalypseLevinEtAl} for every $i \in \{0, \dots, 5\}$ and is given in Table~\ref{tab:RunP2GDR} for some $i \in \{6, \dots, 38\}$.
The only iteration of $\pgdr$ that differs from a $\pgd$ iteration is the fifth one, where $\rank_\Delta X_5 = 1$, $\hat{X}_5^1 = \diag(1, 0, 0)$, and $X_6$ is selected in $\pgd(\hat{X}_5^1; f, \ushort{\alpha}, \bar{\alpha}, \beta, c)$.
The sequences $(f(X_i)-f^*)_{i=0}^{i=38}$, $(\s_f(X_i))_{i=0}^{i=38}$, and $(\norm{X_i-X^*})_{i=0}^{i=38}$ are represented in Figure~\ref{fig:RunP2GDR}.

\begin{table}[h]
\begin{center}
\begin{tabular}{ll}
$i$ & $X_i$\\
\hline
$6$ & $\diag(1.046656000000000, 0, 1.600000000000000)$\\
$11$ & $\diag(1.002866544640000, 0, 1.323933131082407)$\\
$16$ & $\diag(1.000222902511206, 0, 1.324855302786614)$\\
$21$ & $\diag(1.000023110532362, 0, 1.324722970132156)$\\
$26$ & $\diag(1.000001797074997, 0, 1.324717078903522)$\\
$31$ & $\diag(1.000000062106912, 0, 1.324717847681821)$\\
$38$ & $\diag(1.000000002318128, 0, 1.324717955251852)$\\
\hline
\end{tabular}
\end{center}
\caption{Some iterates $X_i$ produced by $\pgdr$ for the problem of Section~\ref{subsec:ExampleApocalypseLevinEtAl}.}
\label{tab:RunP2GDR}
\end{table}

If the behavior of $\pgdr$ on this example seems satisfying, it should however be noted that, if $\Delta < (\frac{3}{5})^{36}$, then $\pgdr$ produces the exact same (finite) sequence of iterates as $\pgd$ because $\rank_\Delta X_{36} = 2$ and $\s_f(X_{37}) \le \varepsilon$. This shows that, in a practical implementation of Algorithm~\ref{algo:IterativeP2GDR} with $\varepsilon > 0$, it is important to choose $\Delta$ in such a way that the algorithm does not stop while it is heading towards an apocalyptic point, which is $\diag(1, 0, 0)$ in this case, in the sense that, if we had continued with $\varepsilon := \Delta := 0$, an apocalypse would have occurred.

\subsection{An apocalypse on $\R_{\le 1}^{2 \times 2}$}
\label{subsec:ExampleApocalypseSmallestSize}
For the function
\begin{equation*}
f : \R^{2 \times 2} \to \R : X \mapsto \frac{X_{1,1}^2+(X_{2,2}-1)^2+(X_{1,2}-X_{2,1})^2}{2},
\end{equation*}
we have $\min_{\R_{\le 1}^{2 \times 2}} f = 0$ and $\argmin_{\R_{\le 1}^{2 \times 2}} f = \diag(0,1)$.
Proposition~\ref{prop:ExampleApocalypseSmallestSize} states that $\pgd$ used with an initial step size for the backtracking procedure smaller than $1$ can follow an apocalypse by trying to minimize $f$ on $\R_{\le 1}^{2 \times 2}$. Before introducing that proposition, we give an intuitive explanation of the result. Given any point $\diag(x_0,0)$ with $x_0 \in (0,\infty)$, $\pgd$ produces a sequence converging to $0_{2 \times 2}$, thereby minimizing the first term of $f$. However, no iteration affects the second term because the search direction $\diag(0,1)$, which would enable the minimization of the second term, is not available until $0_{2 \times 2}$ is reached, which never happens. The third term of $f$ makes its global minimizer on $\R_{\le 1}^{2 \times 2}$ unique without affecting the iterations.

\begin{prop}
\label{prop:ExampleApocalypseSmallestSize}
Let $x_0 \in (0,\infty)$ and $\alpha \in (0,1)$.
With $f$ on $\R_{\le 1}^{2 \times 2}$ as defined above, starting from $X_0 := \diag(x_0,0)$, and using $\bar{\alpha} := \ushort{\alpha} := \alpha$, $\beta \in (0,1)$, and $c \in (0,\frac{1}{2}]$, $\pgd$ produces the sequence $(X_i)_{i \in \N}$ defined by
\begin{equation}
\label{eq:ExampleApocalypseSmallestSize}
X_i := \diag((1-\alpha)^ix_0,0)
\end{equation}
for every $i \in \N$. Moreover, $\s_f(X_i) = (1-\alpha)^ix_0$ for every $i \in \N$. In particular, since $\s_f(0_{2 \times 2}) = \norm{\nabla f(0_{2 \times 2})} = 1$, $(0_{2 \times 2}, (X_i)_{i \in \N}, f)$ is an apocalypse.
\end{prop}

\begin{pf}
The formula \eqref{eq:ExampleApocalypseSmallestSize} holds for $i = 0$. 
Furthermore, for every $X \in \R^{2 \times 2}$,
\begin{equation*}
\nabla f(X) = X - \begin{bmatrix} 0 & X_{2,1} \\ X_{1,2} & 1 \end{bmatrix}.
\end{equation*}
Therefore, for every $i \in \N$,
\begin{align*}
-\nabla f(X_i) = \diag(-(1-\alpha)^ix_0,1),\\
\proj{\tancone{\R_{\le 1}^{2 \times 2}}{X_i}}{-\nabla f(X_i)} = \diag(-(1-\alpha)^ix_0,0),
\end{align*}
and the formula for $\s_f(X_i)$ is valid.
Thus, for every $i \in \N$,
\begin{align*}
X_{i+1} &= X_i + \alpha \proj{\tancone{\R_{\le 1}^{2 \times 2}}{X_i}}{-\nabla f(X_i)},\\
f(X_{i+1}) &\le f(X_i) - c \, \alpha \, \s_f(X_i)^2,
\end{align*}
which shows that the sequence defined by \eqref{eq:ExampleApocalypseSmallestSize} is indeed the one produced by $\pgd$.
The expression for $\s_f(0_{2 \times 2})$ follows from the fact that $-\nabla f(0_{2 \times 2}) = \diag(0,1) \in \R_{\le 1}^{2 \times 2} = \tancone{\R_{\le 1}^{2 \times 2}}{0_{2 \times 2}}$.
\hfill$\qed$
\end{pf}

The next proposition shows that $\pgdr$ escapes the apocalypse due to its rank reduction mechanism. During the first iterations, $\pgdr$ produces the same iterates as $\pgd$. However, when the numerical rank of the iterate becomes smaller than its rank, i.e., when its smallest singular value becomes smaller than or equal to $\Delta$, $\pgdr$ realizes that a stronger decrease of $f$ is obtained by first reducing the rank and then applying an iteration of $\pgd$. As a result, the first term of $f$ is minimized within a finite number of iterations, after which the minimization of the second term can start.

\begin{prop}
\label{prop:ExampleApocalypseSmallestSizeP2GDR}
Consider the same problem as in Proposition~\ref{prop:ExampleApocalypseSmallestSize} with the same parameters and $\Delta \in (0,\infty)$.
Then, $\pgdr$ produces the sequence $(X_i)_{i \in \N}$ defined by
\begin{equation}
\label{eq:ExampleApocalypseSmallestSizeP2GDR}
X_i := \left\{\begin{array}{ll}
\diag((1-\alpha)^ix_0,0) & \text{if } i \le i_\Delta\\
\diag(0,1-(1-\alpha)^{i-i_\Delta}) & \text{if } i > i_\Delta\\
\end{array}\right.
\end{equation}
where $i_\Delta := \max\Big\{\Big\lceil\frac{\ln(\frac{\Delta}{x_0})}{\ln(1-\alpha)}\Big\rceil,0\Big\}$.
In particular, $(X_i)_{i \in \N}$ converges to $\diag(0,1)$ and $\lim_{i \to \infty} \s_f(X_i) = 0$.
\end{prop}

\begin{pf}
The formula \eqref{eq:ExampleApocalypseSmallestSizeP2GDR} is correct for $i = 0$.
If $i_\Delta > 0$, then $(1-\alpha)^ix_0 > \Delta$ for every $i \in \{0, \dots, i_\Delta-1\}$, and \eqref{eq:ExampleApocalypseSmallestSizeP2GDR} thus holds for every $i \in \{1, \dots, i_\Delta\}$ in view of Proposition~\ref{prop:ExampleApocalypseSmallestSize}.
It remains to prove \eqref{eq:ExampleApocalypseSmallestSizeP2GDR} for every integer $i > i_\Delta$.
Let us look at iteration $i_\Delta$. Since $\hat{X}_{i_\Delta}^1 = 0_{2 \times 2}$, $-\nabla f(0_{2 \times 2}) = \diag(0,1) \in \R_{\le 1}^{2 \times 2} = \tancone{\R_{\le 1}^{2 \times 2}}{0_{2 \times 2}}$, $\s_f(0_{2 \times 2}) = 1$, $\hat{X}_{i_\Delta}^1 - \alpha \nabla f(\hat{X}_{i_\Delta}^1) = \diag(0,\alpha)$ and
\begin{equation*}
f(\hat{X}_{i_\Delta}^1)-f(\diag(0,\alpha)) \ge c \, \alpha \, \s_f(\hat{X}_{i_\Delta}^1)^2,
\end{equation*}
we have $\tilde{X}_{i_\Delta}^1 = \diag(0,\alpha)$.
As $\hat{X}_{i_\Delta}^0 = X_{i_\Delta}$, Proposition~\ref{prop:ExampleApocalypseSmallestSize} yields $\tilde{X}_{i_\Delta}^0 = \diag((1-\alpha)^{i_\Delta+1}x_0,0)$. Since
\begin{equation*}
f(\tilde{X}_{i_\Delta}^1)
= \frac{(1-\alpha)^2}{2}
< \frac{(1-\alpha)^{2(i_\Delta+1)}x_0^2+1}{2}
= f(\tilde{X}_{i_\Delta}^0),
\end{equation*}
we have $X_{i_\Delta+1} = \tilde{X}_{i_\Delta}^1$, in agreement with \eqref{eq:ExampleApocalypseSmallestSizeP2GDR}.
Let us now assume that \eqref{eq:ExampleApocalypseSmallestSizeP2GDR} holds for some integer $i > i_\Delta$ and prove that it also holds for $i+1$. As $\hat{X}_i^0 = X_i$, $-\nabla f(X_i) = \diag(0,(1-\alpha)^{i-i_\Delta}) \in \tancone{\R_{\le 1}^{2 \times 2}}{X_i}$, $\s_f(X_i) = (1-\alpha)^{i-i_\Delta}$, $X_i - \alpha \nabla f(X_i) = \diag(0,1-(1-\alpha)^{i+1-i_\Delta}) \in \R_{\le 1}^{2 \times 2}$, and $f(X_i) - f(\diag(0,1-(1-\alpha)^{i+1-i_\Delta})) \ge c \, \alpha \, \s_f(X_i)^2$, we have $\tilde{X}_i^0 = \diag(0,1-(1-\alpha)^{i+1-i_\Delta})$. If $\rank_\Delta X_i = 0$, then $\pgdr$ also considers $\hat{X}_i^1 = 0_{2 \times 2}$ and, from what precedes, $\tilde{X}_i^1 = \diag(0,\alpha)$. Since $f(\tilde{X}_i^0) < f(\tilde{X}_i^1)$, we have $X_{i+1} = \tilde{X}_i^0$, as wished. The other two claims follow.
\hfill$\qed$
\end{pf}

\begin{figure}[h]
\begin{center}
\begin{tikzpicture}
\def\x{1}
\def\a{0.6}
\def\D{0.2}
\pgfmathsetmacro\iD{ceil(ln(\D/\x)/ln(1-\a))}
\draw [->] (0,0) -- (2,0) node [right] {$x$};
\draw [->] (0,0) -- (0,3) node [above] {$y$};
\foreach \i in {0, 1, ..., 4}
{
\draw [dashed] (0,1) circle ({0.4*\i});
}
\foreach \i in {0, 1, ..., \iD}
{
\draw ({\x*(1-\a)^\i},0) node {$\boldsymbol\cdot$};
\draw [->] ({\x*(1-\a)^\i},0) -- ({\x*(1-\a)^(\i+1)},\a);
}
\foreach \i in {1, 2, ..., 3}
{
\draw (0,{1-(1-\a)^\i},0) node {$\boldsymbol\cdot$};
}
\end{tikzpicture}
\end{center}
\caption{Iterates $X_i$ produced by $\pgdr$ for the problem of Section~\ref{subsec:ExampleApocalypseSmallestSize} with $x_0 := 1$, $\alpha := \frac{3}{5}$, and $\Delta := \frac{1}{5}$ in the $xy$-plane of $\diag(x,y)$ matrices. The arrows represent $-\alpha\nabla f(X_i)$.}
\label{fig:ExampleApocalypseSmallestSize}
\end{figure}
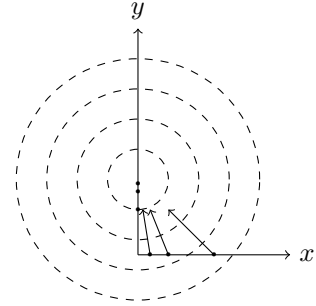

The iterates of $\pgdr$ computed in Proposition~\ref{prop:ExampleApocalypseSmallestSizeP2GDR} are represented in Figure~\ref{fig:ExampleApocalypseSmallestSize}, which summarizes this subsection. As explained, $\pgd$ follows an apocalypse because, at any point $\diag(x,0)$ with $x \in (0,\infty)$, the projection of $-\nabla f$ onto the tangent cone to $\R_{\le 1}^{2 \times 2}$ is parallel to the $x$-axis, and can thus minimize only the first term of $f$. The descent direction $\diag(0,1)$, which enables the minimization of the second term of $f$, becomes accessible only at $\diag(0,0)$.

Although $\pgdr$ avoids the apocalypse for every $\Delta > 0$, it should be noted that, if $\Delta \ge \alpha$, then its rank reduction mechanism makes it apply the $\pgd$ map to $0_{2 \times 2}$ in at least one iteration from iteration $i_\Delta+1$, thereby constructing points that are not used, as shown in the proof of Proposition~\ref{prop:ExampleApocalypseSmallestSizeP2GDR}. For those iterations, $\pgdr$ therefore produces the same iterates as $\pgd$ at a higher computational cost.

We close this section by discussing how Algorithm~\ref{algo:IterativeP2GDR} with $\varepsilon > 0$ behaves on this problem. If $\Delta := 0$, it returns the sequence $(X_i)_{i=0}^{i=i_\varepsilon}$ defined by \eqref{eq:ExampleApocalypseSmallestSize}, where $i_\varepsilon := \max\Big\{\Big\lceil\frac{\ln(\frac{\varepsilon}{x_0})}{\ln(1-\alpha)}\Big\rceil,0\Big\}$. Thus, in view of \eqref{eq:ExampleApocalypseSmallestSizeP2GDR}, for Algorithm~\ref{algo:IterativeP2GDR} to avoid stopping while it is heading towards the apocalyptic point, we must have $i_\Delta < i_\varepsilon$, i.e., $\Delta \ge (1-\alpha)^{i_\varepsilon-1}x_0$.

\section{A potential side effect of the rank reduction mechanism}
\label{sec:PotentialSideEffectRankReductionMechanism}
In this section, we report our observation that, for some instances of \eqref{eq:LowRankOpti}, $\pgd$ and $\pgdr$ may converge to different stationary points with different costs, especially if $\Delta$ is large. In Section~\ref{subsec:P2GD<P2GDR}, we present an example where the stationary point to which $\pgd$ converges has a lower cost than the one to which $\pgdr$ converges. The converse situation happens in the example given in Section~\ref{subsec:P2GD>P2GDR}. In both examples, the initial iterate is $X_0 := \diag(1,0)$, $\pgdr$ is used with $\Delta := 1$, and we have $0 = \rank_\Delta X_0 < \rank X_0 = 1$. Moreover, for the two proposed cost functions, applying an iteration of $\pgd$ to $\hat{X}_0^1 := 0_{2 \times 2}$ decreases the cost more than to $\hat{X}_0^0 := X_0$, and $\pgdr$ therefore computes $X_1$ by applying a $\pgd$ iteration to $0_{2 \times 2}$. In this case, using $\pgdr$ with $\Delta := 1$ amounts to starting $\pgd$ from another initial iterate, and hence it is not surprising that another output is produced.

\subsection{$\pgd$ may find a better stationary point than $\pgdr$}
\label{subsec:P2GD<P2GDR}
Define the function $f : \R^{2 \times 2} \to \R$ by
\begin{equation*}
f(X) := \frac{(X_{1,1}-4)^2+3(X_{2,2}-2)^2+(X_{1,2}-X_{2,1})^2}{2}.
\end{equation*}
Then,
\begin{align*}
\min_{\R_{\le 1}^{2 \times 2}} f = 0,&&
\argmin_{\R_{\le 1}^{2 \times 2}} f = \left\{\begin{bmatrix} 4 & \pm 2\sqrt{2} \\ \pm 2\sqrt{2} & 2 \end{bmatrix}\right\}.
\end{align*}
The next two propositions show that, for some input parameters, $\pgd$ converges to a stationary point with lower cost than the one to which $\pgdr$ converges.
Starting from $\diag(1,0)$, neither of the two algorithms converges to one of the two global minimizers of $f$ on $\R_{\le 1}^{2 \times 2}$. On the one hand, $\pgd$ produces a sequence converging to $\diag(4,0)$, thereby minimizing the first term of $f$ and achieving a cost of $6$. On the other hand, at the first iteration, $\pgdr$ prefers to apply the $\pgd$ map to $0_{2 \times 2}$ because this yields a stronger decrease of $f$ thanks to the factor $3$ in the second term. After that first iteration, $\pgdr$ produces the same iterates as $\pgd$ would and constructs a sequence converging to $\diag(0,2)$, thus minimizing the second term and achieving a cost of $8$. The third term of $f$ makes its set of minimizers on $\R_{\le 1}^{2 \times 2}$ finite without affecting the iterations of $\pgd$ and $\pgdr$.

\begin{prop}
\label{prop:P2GD<P2GDR}
Let $\alpha := \frac{1}{4}$.
With $f$ on $\R_{\le 1}^{2 \times 2}$ as defined above, starting from $X_0 := \diag(1,0)$, and using $\bar{\alpha} := \ushort{\alpha} := \alpha$, $\beta \in (0,1)$, and $c \in (0,\frac{5}{8}]$, $\pgd$ produces the sequence $(X_i)_{i \in \N}$ defined by $X_i := \diag(4-3(1-\alpha)^i,0)$ for every $i \in \N$.
In particular, $(X_i)_{i \in \N}$ converges to $\diag(4,0)$, $\lim_{i \to \infty} \s_f(X_i) = 0$, and $\s_f(\diag(4,0)) = 0$.
\end{prop}

\begin{pf}
The formula holds for $i = 0$. Let us prove that, if it holds for some $i \in \N$, it also holds for $i+1$.
For every $X \in \R^{2 \times 2}$,
\begin{equation*}
\nabla f(X) = \begin{bmatrix} X_{1,1}-4 & X_{1,2}-X_{2,1} \\ X_{2,1}-X_{1,2} & 3(X_{2,2}-2) \end{bmatrix}.
\end{equation*}
Thus,
\begin{align*}
-\nabla f(X_i) = \diag(3(1-\alpha)^i,6),\\
\proj{\tancone{\R_{\le 1}^{2 \times 2}}{X_i}}{-\nabla f(X_i)} = \diag(3(1-\alpha)^i,0),
\end{align*}
and $\s_f(X_i) = 3(1-\alpha)^i$. Since
\begin{equation*}
X_i + \alpha \diag(3(1-\alpha)^i,0) = \diag(4-3(1-\alpha)^{i+1},0) \in \R_{\le 1}^{2 \times 2}
\end{equation*}
and
\begin{equation*}
f(X_i) - f(\diag(4-3(1-\alpha)^{i+1},0)) \ge c \, \alpha \, \s_f(X_i)^2,
\end{equation*}
$X_{i+1}$ has the required form. The other claims follow.$\qed$
\end{pf}

\begin{prop}
Consider the same problem as in Proposition~\ref{prop:P2GD<P2GDR} with the same parameters and $\Delta := 1$.
Then, $\pgdr$ produces the sequence $(X_i)_{i \in \N}$ defined by $X_i := \diag(0,2-2(1-3\alpha)^i)$ for every $i \in \N \setminus \{0\}$.
In particular, $(X_i)_{i \in \N}$ converges to $\diag(0,2)$, $\lim_{i \to \infty} \s_f(X_i) = 0$, and $\s_f(\diag(0,2)) = 0$.
\end{prop}

\begin{pf}
Let us first prove the formula for $i = 1$. By the proof of the preceding proposition, as $\hat{X}_0^0 = X_0$, we have $\tilde{X}_0^0 = \diag(1+3\alpha,0)$. Since $\hat{X}_0^1 = 0_{2 \times 2}$, $-\nabla f(0_{2 \times 2}) = \diag(4,6)$, $\proj{\tancone{\R_{\le 1}^{2 \times 2}}{0_{2 \times 2}}}{-\nabla f(0_{2 \times 2})} = \proj{\R_{\le 1}^{2 \times 2}}{-\nabla f(0_{2 \times 2})} = \diag(0,6)$, $\s_f(0_{2 \times 2}) = 6$, $\hat{X}_0^1 + \alpha \diag(0,6) = \diag(0,6\alpha)$, and
\begin{equation*}
f(\hat{X}_0^1)-f(\diag(0,6\alpha)) \ge c \, \alpha \, \s_f(\hat{X}_0^1)^2,
\end{equation*}
we have $\tilde{X}_0^1 = \diag(0,6\alpha)$.
As
\begin{equation*}
f(\tilde{X}_0^1)
= \frac{16+12(1-3\alpha)^2}{2}
< \frac{9(1-\alpha)^2+12}{2}
= f(\tilde{X}_0^0),
\end{equation*}
we have $X_1 = \tilde{X}_0^1$, in agreement with the formula.
Let us now assume that the formula holds for some $i \in \N \setminus \{0\}$ and prove that it also holds for $i+1$. Observe that $\rank_\Delta X_i = 1$. Since $-\nabla f(X_i) = \diag(4,6(1-3\alpha)^i)$, $\proj{\tancone{\R_{\le 1}^{2 \times 2}}{X_i}}{-\nabla f(X_i)} = \diag(0,6(1-3\alpha)^i)$, $\s_f(X_i) = 6(1-3\alpha)^i$, $X_i + \alpha \diag(0,6(1-3\alpha)^i) = \diag(0,2-2(1-3\alpha)^{i+1}) \in \R_{\le 1}^{2 \times 2}$, and
\begin{equation*}
f(X_i) - f(\diag(0,2-2(1-3\alpha)^{i+1})) \ge c \, \alpha \, \s_f(X_i)^2,
\end{equation*}
the formula is valid for $i+1$. The other claims follow.$\qed$
\end{pf}

\subsection{$\pgdr$ may find a better stationary point than $\pgd$}
\label{subsec:P2GD>P2GDR}
Define the function $f : \R^{2 \times 2} \to \R$ by
\begin{equation*}
f(X) := \frac{(X_{1,1}-2)^2+(X_{2,2}-3)^2+(X_{1,2}-X_{2,1})^2}{2}.
\end{equation*}
Then,
\begin{align*}
\min_{\R_{\le 1}^{2 \times 2}} f = 0,&&
\argmin_{\R_{\le 1}^{2 \times 2}} f = \left\{\begin{bmatrix} 2 & \pm \sqrt{6} \\ \pm \sqrt{6} & 3 \end{bmatrix}\right\}.
\end{align*}
The next two propositions show that, for some input parameters, $\pgd$ and $\pgdr$ achieve the costs of $\frac{9}{2}$ and $2$, respectively.

\begin{prop}
\label{prop:P2GD>P2GDR}
Let $\alpha \in (\frac{1}{3},1)$.
With $f$ on $\R_{\le 1}^{2 \times 2}$ as defined above, starting from $X_0 := \diag(1,0)$, and using $\bar{\alpha} := \ushort{\alpha} := \alpha$, $\beta \in (0,1)$, and $c \in (0,\frac{1}{2}]$, $\pgd$ produces the sequence $(X_i)_{i \in \N}$ defined by $X_i := \diag(2-(1-\alpha)^i,0)$ for every $i \in \N$.
In particular, $(X_i)_{i \in \N}$ converges to $\diag(2,0)$, $\lim_{i \to \infty} \s_f(X_i) = 0$, and $\s_f(\diag(2,0)) = 0$.
\end{prop}

\begin{pf}
The formula holds for $i = 0$. Let us prove that, if it holds for some $i \in \N$, it also holds for $i+1$.
For every $X \in \R^{2 \times 2}$,
\begin{equation*}
\nabla f(X) = X - \begin{bmatrix} 2 & X_{2,1} \\ X_{1,2} & 3 \end{bmatrix}.
\end{equation*}
Thus,
\begin{align*}
-\nabla f(X_i) = \diag((1-\alpha)^i,3),\\
\proj{\tancone{\R_{\le 1}^{2 \times 2}}{X_i}}{-\nabla f(X_i)} = \diag((1-\alpha)^i,0),
\end{align*}
and $\s_f(X_i) = (1-\alpha)^i$. Since
\begin{equation*}
X_i + \alpha \diag((1-\alpha)^i,0) = \diag(2-(1-\alpha)^{i+1},0) \in \R_{\le 1}^{2 \times 2}
\end{equation*}
and
\begin{equation*}
f(X_i) - f(\diag(2-(1-\alpha)^{i+1},0)) \ge c \, \alpha \, \s_f(X_i)^2,
\end{equation*}
$X_{i+1}$ has the required form. The other claims follow.$\qed$
\end{pf}

\begin{prop}
Consider the same problem as in Proposition~\ref{prop:P2GD>P2GDR} with the same parameters and $\Delta := 1$.
Then, $\pgdr$ produces the sequence $(X_i)_{i \in \N}$ defined by $X_i := \diag(0,3-3(1-\alpha)^i)$ for every $i \in \N \setminus \{0\}$.
In particular, $(X_i)_{i \in \N}$ converges to $\diag(0,3)$, $\lim_{i \to \infty} \s_f(X_i) = 0$, and $\s_f(\diag(0,3)) = 0$.
\end{prop}

\begin{pf}
Let us first prove the formula for $i = 1$. By the proof of the preceding proposition, as $\hat{X}_0^0 = X_0$, we have $\tilde{X}_0^0 = \diag(1+\alpha,0)$. Since $\hat{X}_0^1 = 0_{2 \times 2}$, $-\nabla f(0_{2 \times 2}) = \diag(2,3)$, $\proj{\tancone{\R_{\le 1}^{2 \times 2}}{0_{2 \times 2}}}{-\nabla f(0_{2 \times 2})} = \proj{\R_{\le 1}^{2 \times 2}}{-\nabla f(0_{2 \times 2})} = \diag(0,3)$, $\s_f(0_{2 \times 2}) = 3$, $\hat{X}_0^1 + \alpha \diag(0,3) = \diag(0,3\alpha)$, and
\begin{equation*}
f(\hat{X}_0^1)-f(\diag(0,3\alpha))
\ge c \, \alpha \, \s_f(\hat{X}_0^1)^2,
\end{equation*}
we have $\tilde{X}_0^1 = \diag(0,3\alpha)$.
As
\begin{equation*}
f(\tilde{X}_0^1)
= \frac{4+9(1-\alpha)^2}{2}
< \frac{(1-\alpha)^2+9}{2}
= f(\tilde{X}_0^0),
\end{equation*}
we have $X_1 = \tilde{X}_0^1$, in agreement with the formula.
Let us now assume that the formula holds for some $i \in \N \setminus \{0\}$ and prove that it also holds for $i+1$. Observe that $\rank_\Delta X_i = 1$. Since $-\nabla f(X_i) = \diag(2,3(1-\alpha)^i)$, $\proj{\tancone{\R_{\le 1}^{2 \times 2}}{X_i}}{-\nabla f(X_i)} = \diag(0,3(1-\alpha)^i)$, $\s_f(X_i) = 3(1-\alpha)^i$, $X_i + \alpha \diag(0,3(1-\alpha)^i) = \diag(0,3-3(1-\alpha)^{i+1}) \in \R_{\le 1}^{2 \times 2}$, and
\begin{equation*}
f(X_i) - f(\diag(0,3-3(1-\alpha)^{i+1})) \ge c \, \alpha \, \s_f(X_i)^2,
\end{equation*}
the formula is valid for $i+1$. The other claims follow.$\qed$
\end{pf}

\section{Conclusion}
\label{sec:Conclusion}
This paper compares $\pgd$ and $\pgdr$ on synthetic instances of \eqref{eq:LowRankOpti}. The simplicity of those instances enables both analytical and empirical investigations. This allows us to observe two behaviors:
\begin{enumerate}
\item $\pgd$ and $\pgdr$ respectively following and escaping apocalypses (Section~\ref{sec:TwoExamplesApocalypses});
\item $\pgd$ and $\pgdr$ converging to different stationary points with different cost values (Section~\ref{sec:PotentialSideEffectRankReductionMechanism}).
\end{enumerate}

Concerning $\pgdr$, i.e., Algorithm~\ref{algo:IterativeP2GDR} with $\Delta > 0$, we also observe that, if $\varepsilon > 0$, which is always the case in a practical implementation, then the choice of $\Delta$ can be significant. Indeed, when $\varepsilon := 0$, the convergence analysis of Algorithm~\ref{algo:IterativeP2GDR} given in \cite{OlikierGallivanAbsil2022} holds for every $\Delta > 0$. However, as remarked in Section~\ref{sec:TwoExamplesApocalypses}, if $\varepsilon > 0$, then $\Delta$ must be chosen large enough to let the rank reduction mechanism prevent the algorithm from stopping while heading towards an apocalyptic point. On the other hand, choosing $\Delta$ too large can make the rank reduction mechanism work inefficiently in the sense that, for some iterations, $\pgdr$ produces the same iterates as $\pgd$ at a higher computational cost. In any case, it is good practice, in order to avoid following an apocalypse, to apply a rank reduction to the last iterate and to look at the effect of a $\pgd$ iteration on the obtained point.
Besides these apocalypse-related considerations, the larger $\Delta$ is, the more side effects of the rank reduction mechanism are likely to arise, as noticed in Section~\ref{sec:PotentialSideEffectRankReductionMechanism}.

We close this paper with two open questions regarding $\pgd$ and $\pgdr$.
\begin{enumerate}
\item As pointed out at the end of Section~\ref{sec:NotationPreliminaries}, it is not known whether there exists an instance of \eqref{eq:LowRankOpti} for which $\pgd$ produces a sequence with the following two properties: it does not diverge to infinity and the stationarity measure $\s_f$ does not go to zero along any convergent subsequence. For such an instance, there would exist $\varepsilon > 0$ such that Algorithm~\ref{algo:IterativeP2GDR} with $\Delta := 0$ does not terminate.
\item Is there an instance of \eqref{eq:LowRankOpti} for which $\pgd$ converges to a nonstationary point having a lower cost than the stationary point to which $\pgdr$ converges?
\end{enumerate}

\bibliography{golikier_bib}

\end{document}